\documentclass[12pt,a4paper,reqno]{amsart}

\usepackage{amsmath}
\usepackage{amssymb}
\usepackage{amsfonts}

\advance\textheight by \topskip
\def\O{\mathcal{O}}
\allowdisplaybreaks
\newcommand{\qf}[2]{(#1;q)_{#2}}
\newtheorem{theorem}{Theorem}

\newtheorem{lemma}{Lemma}

\def\la{\lambda}
\def\CO{\mathbb{C}}
\def\Prob{\mathbb{P}}

\def\bb#1{[\![#1]\!]}

\title[ Value and position of large 
left--to--right maxima]
{Combinatorics of geometrically distributed random
variables:\\  
Value and position of large 
left--to--right maxima}

\author{Helmut Prodinger}

\address{ Helmut Prodinger,
The John Knopfmacher Centre for Applicable Analysis and Number Theory,
 Department of Mathematics,
University of the Witwatersrand, P.~O. Wits, 
2050 Johannesburg, South Africa, email:
{\tt helmut@gauss.cam.wits.ac.za},\newline
homepage:
{\tt http://www.wits.ac.za/helmut/index.htm}.
}

\date{\today}

\begin{document}

\begin{abstract}
For words of length $n$, generated by independent geometric
random variables, we consider the average value and the
average position of the $r$th left--to--right maximum
{\it counted from the right,} for
fixed $r$ and $n\to \infty$. This complements previous research
\cite{KnPr00} where the analogous questions were considered
for the $r$th left--to--right maximum
{\it counted from the left.} 
\end{abstract}

\maketitle

\section{Introduction}
\label{sec:intro}

In this paper we study sequences  with respect
to the {\em left--to--right maxima\/} of the elements of
the sequence. 

This paper is a companion paper to \cite{KnPr00}, and these
two papers should be read together.

In \cite{Prodinger96c} the number
of left--to--right maxima was investigated
in the model of {\em words\/} (strings) $a_1\dots a_n$, 
where the letters $a_i\in\mathbb N$ are independently
generated according to the geometric distribution with
\begin{equation}\label{geom}
\Prob{\{X=k\}}=pq^{k-1}, \qquad\text{where $p+q=1$.}
\end{equation}
(An element is a left--to--right maximum if it is 
larger than all the elements to the left. The version
where `larger' is replaced by `larger or equal' would
be also of interest but is not considered here.)

In  \cite{KnPr00} we considered  the two 
parameters `value' and `position' of the 
$r$th left--to--right maximum
for geometric random variables and 
obtained the 
asymptotic
formul{\ae}
$\frac rp$ and 
$\frac{1}{(r-1)!}\big(\frac pq\log_Qn\big)^{r-1}$.

These results cover the instances of small
left--to--right maxima; therefore we investigate
here the instance of {\it large\/} left--to--right maxima.
We say that the last left--to--right maximum is the
instance $r=1$, the previous left--to--right maximum
is the instance $r=2$, and so on.   

As can be seen in the forthcoming book \cite{FlSe05},
in the limit $q\to1$ the model becomes the model
of {\it random permutations.} The concept of `value' does
not carry over, but the  concept of `position' does.

The largest element in a random permutation of $n$ elements
(which is $n$) is expected to be in the middle.
The second largest is then expected to be in the middle
of the first half, and so on. In general, the average
position should be about $\frac{n}{2^r}$. 

We find it useful  to use the additional
notations $Q:=q^{-1}$,  $L:=\log Q$, 
and
$\bb{i}:=1-(1-q^i)z$.

In \cite{Prodinger96c}, compare also
\cite{SzRe90}, the  average value of the {\it height\/}
appeared; it is asymptotic to $\log_Qn$, which is
also the instance $r=1$ of the average value
of left--to--right maxima, counted from the right. 
Intuitively, for fixed $r$, the same asymptotic
formula is to be expected, the depency of $r$    only
appearing in the lower order terms.

It should be noted that not all 
random strings of length $n$
{\em have\/} $r$
left--to--right maxima.
However, as explained already in \cite{KnPr00},
this happens with low probability and can be ignored.

\section{The average value of the $r$th left--to--right
maximum from the right}

The following expression for  
the average value of the $r$th left--to--right
maximum from the right is obtained as in
\cite{KnPr00}; just observe that we assume that it
is $h$, and that $r-1$ more left--to--right
maxima follow to the right, whereas  to the left
everything is smaller than $h$:

\begin{align*}
 E^{(r)}_{n}&=[z^n]
\sum_{1\le h< i_1< \dots <i_{r-1} }
h\frac{1}{\bb{h-1}}\frac{zpq^{h-1}}{\bb{h}}
\frac{zpq^{i_{1}-1}}{\bb{ i_{1}}}\dots
\frac{zpq^{i_{r-1}-1}}{\bb {i_{r-1}}}
\\
&=
\Big(\frac{p}{q}\Big)^{r}[z^n]
\sum_{1\le h< i_1< \dots <i_{r-1} }h
\frac{1}{\bb{h-1}}\frac{zq^h}{\bb{h}}
\frac{q^{i_{1}}z}{\bb{ i_{1}}}\dots
\frac{q^{i_{r-1}}z}{\bb {i_{r-1}}}
\\
&=(-1)^r
\Big(\frac{p}{q}\Big)^{r}(-1)^n[w^n](1-w)^{n}\times
\\&\qquad\qquad
\times
\sum_{1\le h< i_1< \dots <i_{r-1} }\frac{1}{1-wq^{h-1}}
\frac{wq^{h}h}{1-q^hw}
\frac{q^{i_{1}}w\dots q^{i_{r-1}}w}{
(1-q^{i_{1}}w)\dots(1-q^{i_{r-1}}w)}
\\
&=(-1)^r
\Big(\frac{p}{q}\Big)^{r}
\sum_{k=r}^{n}\binom{n}{k}(-1)^k\times\\&\qquad\qquad
\times[w^k]
\sum_{1\le h< i_1< \dots <i_{r-1} }\frac{hwq^{h}}{(1-wq^{h-1})(1-wq^{h})}
\frac{q^{i_{1}}w\dots q^{i_{r-1}}w}{
(1-q^{i_{1}}w)\dots(1-q^{i_{r-1}}w)}
.\\
\end{align*}

For the computation of the inner sum we need the following
formula, which can be computed by elementary manipulations.
For $n\ge m$ we have
\begin{align*}
[w^n]\sum_{h\ge1}\frac{hwq^h}{(1-wq^{h-1})(1-wq^{h})}(wq^h)^m
=\frac qp\frac{q^m-q^n}{(1-q^n)^2}.
\end{align*}

Using this and the principle that

\begin{equation*}
[w^n]\sum_ia_if(b_iw)=\sum_ia_ib_i^n\cdot [w^n]f(w)
\end{equation*} 
we find that

\begin{multline*}
[w^k]
\sum_{1\le h< i_1< \dots <i_{r-1} }\frac{hwq^{h}}{(1-wq^{h-1})(1-wq^{h})}
\frac{q^{i_{1}}w\dots q^{i_{r-1}}w}{
(1-q^{i_{1}}w)\dots(1-q^{i_{r-1}}w)}\\
=\frac{q}{p}\frac{1}{(1-q^k)^2}\sum_{0<l_1<\dots<l_{r-1}<k}
\frac{1}{Q^{l_1}-1}\dots\frac{1}{Q^{l_{r-1}}-1}
\left(q^{l_{r-1}}-q^k\right).\\
\end{multline*}

It is worthwhile to note the instance $r=1$, which is the {\it maximum:}

\begin{align*}
 E^{(1)}_{n}
&=
\sum_{k=1}^{n}\binom{n}{k}(-1)^{k-1}
\frac{1}{1-q^k}.
\end{align*}

This checks with previous results, compare \cite{Prodinger96c, SzRe90}, 
but it seems to be a folklore result. 

Now define the function

\begin{equation*}
\phi_r(k):=\sum_{0<l_1<\dots<l_{r-1}<k}
\frac{1}{Q^{l_1}-1}\dots\frac{1}{Q^{l_{r-1}}-1}
(q^{l_{r-1}}-q^k).
\end{equation*}

With this notation, we have
\begin{align*}
 E^{(r)}_{n}
&=
(-1)^r
\Big(\frac{p}{q}\Big)^{r-1}\sum_{k=r}^{n}\binom{n}{k}(-1)^{k}
\frac{\phi_r(k)}{(1-q^k)^2}.
\end{align*}

Since
our quantities come out as alternating sums, 
the appropriate treatment of them is {\it Rice's method\/}
which is surveyed in \cite{FlSe95}; the key point
is the 
following Lemma.
\begin{lemma} 
Let $\mathcal C$ be a curve surrounding the points $1,2,\dots,n$ 
in the complex plane and let $f(z)$ be analytic inside $\mathcal C$. 
Then
$$
\sum_{k=1}^n \binom nk \, {(-1)}^k f(k)=
-\frac 1{2\pi i} \int_{\mathcal C} [n;z] f(z) dz, 
$$
where
$$
[n;z]=\frac{(-1)^{n-1} n!}{z(z-1)\dots(z-n)}=
\frac{\Gamma(n+1)\Gamma(-z)}{\Gamma(n+1-z)}.
$$

\end{lemma}

Extending the contour of integration it turns out that 
under suitable growth conditions on $f(z)$ (compare \cite{FlSe95})
the asymptotic expansion
of the alternating sum is given by
$$ 
\sum \text {Res} \big( [n;z] f(z)\big)+\text{smaller order terms}
$$
where the sum is taken over all poles $z_0$
different from $1,\dots,n$. 
Poles that lie more to the left lead to smaller terms in the
asymptotic expansion.

The range $1,\dots,n$ for the summation is not sacred; 
if we sum, for example, over $k=2,\dots,n$, the contour must
encircle $2,\dots,n$, etc. \qed

\medskip

Therefore a key issue is to extend sequences like
$\phi_r(k)$ to the complex plane, in other words,
give a meaningful definition of  
$\phi_r(z)$, for $z\in\CO$. 

Note that
\begin{align*}
\sum_{0<i_1<\dots<i_{s}<j}\frac{1}{Q^{i_1}-1}\dots
\frac{1}{Q^{i_s}-1}&=[x^{s}]
\prod_{i=1}^{j-1}\left(1+\frac{x}{Q^i-1}\right)\\
&=[x^{s}]
\prod_{i=1}^{j-1}\frac{Q^i-1+x}{Q^i-1}\\
&=[x^{s}]
\prod_{i=1}^{j-1}\frac{1-(1-x)q^i}{1-q^i}\\
&=[x^{s}]
\frac{\qf{(1-x)q}{j-1}}{\qf q{j-1}}.\\
\end{align*}

In this computation we use the standard notation for $q$--factorials,
\begin{equation*}   
\qf an:=(1-a)(1-aq)\dots (1-aq^{n-1}),\qquad
\qf a\infty:=(1-a)(1-aq)\dots,
\end{equation*} 
see e.~g. \cite{Andrews76}. In the following we need
Cauchy's theorem ($q$--binomial theorem), viz.
\begin{equation*}   
\sum_{n\ge0}\frac{\qf an}{\qf qn}t^n=\frac{\qf {at}\infty}
{\qf {t}\infty}.
\end{equation*}

Then
\begin{align*}
\phi_r(k)&=[x^{r-2}]\sum_{1\le l \le k  }
\frac{\qf{(1-x)q}{l-1}}{\qf q{l-1}}
\frac{q^l-q^k}{Q^l-1}\\
&=[x^{r-2}]\sum_{1\le l \le k  }
\frac{\qf{(1-x)q}{l-1}}{\qf q{l}}
(q^{2l}-q^{k+l})\\
&=[x^{r-2}][t^k]\frac{1}{1-t}\sum_{l\ge1  }
\frac{\qf{(1-x)q}{l-1}}{\qf q{l}}
(q^{2l}-q^{k+l})t^l\\
&=[x^{r-2}][t^k]\frac{1}{1-t}\frac{1}{x}\bigg(-(1-q^k)+
\sum_{l\ge0  }
\frac{\qf{1-x}{l}}{\qf q{l}}
(q^{2l}-q^{k+l})t^l\bigg)\\
&=[x^{r-1}][t^k]\frac{1}{1-t}\bigg(-(1-q^k)
+\frac{\qf{(1-x)q^2t)}{\infty}}{\qf{q^2t}{\infty}}-
q^k\frac{\qf{(1-x)qt}{\infty}}{\qf{qt}{\infty}}
\bigg)\\
&=-[x^{r-1}][t^k]\frac{1}{1-t}\bigg((1-q^k)+
(q^k-1)\frac{\qf{(1-x)qt}{\infty}}{\qf{qt}{\infty}}+
\frac{\qf{(1-x)q^2t}{\infty}}{\qf{qt}{\infty}}xqt
\bigg)\\
&=-\delta_{r1}(1-q^k)-[x^{r-1}][t^k]\bigg(
(q^k-1)\frac{\qf{(1-x)qt}{\infty}}{\qf{t}{\infty}}+
\frac{\qf{(1-x)q^2t}{\infty}}{\qf{t}{\infty}}xqt
\bigg)\\
&=-\delta_{r1}(1-q^k)-[x^{r-1}]\bigg(
(q^k-1)\frac{\qf{(1-x)q}k}{\qf qk}
+xq\frac{\qf{(1-x)q^2}{k-1}}{\qf q{k-1}}
\bigg)\\
&=-\delta_{r1}(1-q^k)-[x^{r-1}]\bigg(
(q^k-1)\frac{\qf{(1-x)q}k}{\qf qk}
+\frac{xq(1-q^k)}{1-(1-x)q}\frac{\qf{(1-x)q^2}{k}}{\qf q{k}}
\bigg).\\
\end{align*}
Hence
\begin{align*}
\frac{\phi_r(k)}{1-q^k}&=-\delta_{r1}+[x^{r-1}]\bigg(
\frac{\qf{(1-x)q}k}{\qf qk}
-\frac{xq}{1-(1-x)q}\frac{\qf{(1-x)q^2}{k}}{\qf q{k}}
\bigg)\\
\end{align*}
and
\begin{align*}
\frac{\phi_r(k)}{1-q^k}\bigg|_{k=0}&=-\delta_{r1}+[x^{r-1}]\bigg(
1-\frac{xq}{1-(1-x)q}
\bigg)\\
&=-[x^{r-1}]\frac{xq}{1-(1-x)q}\\
&=-\frac{q}{p}[x^{r-2}]\frac{1}{1+\frac{xq}{p}}\\
&= (-1)^{r-1}\Big(\frac{q}{p}\Big)^{r-1}.\\
\end{align*}

To get lower order terms, we first note that
\begin{equation*}   
\frac{\qf{(1-x)q}z}{\qf qz}\sim 1+L
\sum_{j\ge1}\frac{xq^j}{(1-q^j)(1-(1-x)q^j)}\cdot z 
\qquad (z\to0)
\end{equation*}
and (for $r\ge2$)
\begin{equation*}   
[x^{r-1}]L
\sum_{j\ge1}\frac{xq^j}{(1-q^j)(1-(1-x)q^j)}\cdot z 
=L(-1)^r\sum_{j\ge1}\frac{q^{j(r-1)}}{(1-q^j)^r}\cdot z. 
\end{equation*}
Also note that
\begin{align*}
\frac{xq}{1-(1-x)q}&\frac{\qf{(1-x)q^2}{z}}{\qf q{z}}
\\&\sim \frac{xq}{1-(1-x)q}+
\frac{xq}{1-(1-x)q}L
\sum_{j\ge1}\frac{xq^{j+1}}{(1-q^j)(1-(1-x)q^{j+1})}\cdot z \\
&=\frac{xq}{1-(1-x)q}+
{xqL}
\sum_{j\ge1}
\frac{q^j}{(1-q^j)^2}
\left[-
\frac{1}{1+\frac{qx}{p}}+
\frac{1}{1+\frac{q^{j+1}x}{1-q^{j+1}}}
\right]
\cdot z \\
\end{align*}
and (for $r\ge2$)
\begin{multline*}   
[x^{r-1}]{xqL}
\sum_{j\ge1}
\frac{q^j}{(1-q^j)^2}
\left[-
\frac{1}{1+\frac{qx}{p}}+
\frac{1}{1+\frac{q^{j+1}x}{1-q^{j+1}}}
\right]
\cdot z
\\=qL(-1)^{r-1}
\sum_{j\ge1}
\frac{q^j}{(1-q^j)^2}
\left[\Big(\frac{q}{p}\Big)^{r-2}-
\Big(\frac{q^{j+1}}{1-q^{j+1}}\Big)^{r-2}
\right]
\cdot z.
\end{multline*}
Hence

\begin{align*}
\frac{\phi_r(z)}{1-q^z}
&\sim (-1)^{r-1}\Big(\frac{q}{p}\Big)^{r-1}+\beta_r\cdot z
\qquad (z\to 0),
\\
\end{align*}
with 
\begin{align*}
\beta_r&=L(-1)^r\sum_{j\ge1}\frac{q^{j(r-1)}}{(1-q^j)^r}
+qL(-1)^{r}
\sum_{j\ge1}
\frac{q^j}{(1-q^j)^2}
\left[\Big(\frac{q}{p}\Big)^{r-2}-
\Big(\frac{q^{j+1}}{1-q^{j+1}}\Big)^{r-2}
\right].
\end{align*}

Computing the residue of 
\begin{equation*}   
\frac{1}{z}\frac{1}{1-q^z}\frac{\phi_r(z)}{1-q^z}
\end{equation*}
at $z=0$ eventually leads to 

\begin{theorem}
The average value of the $r$th left--to--right
maximum from the right is given by
\begin{equation*}   
E^{(r)}_n=\log_Qn+\sigma_r+\delta_r(\log_Qn)+o(1),
\end{equation*}
where $\delta_r(x)$ is a small periodic function of\/
$\log_Qn$ (originating from the poles at
$2k\pi i/\log Q$, $k\in\mathbb{Z}\setminus\{0\}$). 
The constant $\sigma_r$  is given by
\begin{equation*}   
\sigma_r=\frac{\gamma}{L}+\frac{1}{2}-
p^{r-1}
\sum_{j\ge1}\frac{q^{(j-1)(r-1)}}{(1-q^j)^r}
-{p}
\sum_{j\ge1}
\frac{q^{j}}{(1-q^j)^2}
+p^{r-1}
\sum_{j\ge1}
\frac{q^{j(r-1)}}{(1-q^j)^2(1-q^{j+1})^{r-2}}.
\end{equation*}
This holds for $r=1$ as well, when the 3 sums combined
evaluate to 0. 
\end{theorem}

For interest, let us note the following limit:
\begin{equation*}
\lim_{q\to1}(1-q)\sigma_r=
\sum_{k=3}^{r-3}(r-1-k)\zeta(k)+(r-1)\zeta(2)+\zeta(r-2)-\zeta(r)
-\frac{r(r-3)}{2}.
\end{equation*}

\section{The average position of the $r$th left--to--right
maximum from the right}

We consider the generating function

\begin{equation*}
\sum_{h\ge1}\frac{zpq^{h-1}}{(1-zv(1-q^{h-1}))\bb h}
\sum_{h<i_1<\dots <i_{r-1}}\frac{zpq^{i_{1}-1}}{\bb{i_1}}\dots
\frac{zpq^{i_{r-1}-1}}{\bb{i_{r-1}}},
\end{equation*}
where $v$ marks the elements to the left of the  
$r$th left--to--right
maximum from the right. The actual position
is then one more than the number of $v$'s, but
that is overshadowed by the error term anyway. 
To get the average, we must differentiate this
with respect to $v$, and then replace $v$    by~1.
The result is

\begin{align*}
&\sum_{h\ge1}\frac{zpq^{h-1}z(1-q^{h-1})}{\bb{h-1}^2\bb h}
\sum_{h<i_1<\dots <i_{r-1}}\frac{zpq^{i_{1}-1}}{\bb{i_1}}\dots
\frac{zpq^{i_{r-1}-1}}{\bb{i_{r-1}}}\\*
&=\Big(\frac pq\Big)^r
\sum_{h\ge1}\frac{zq^{h}z(1-q^{h-1})}{\bb{h-1}^2\bb h}
\sum_{h<i_1<\dots <i_{r-1}}\frac{zq^{i_{1}}}{\bb{i_1}}\dots
\frac{zq^{i_{r-1}}}{\bb{i_{r-1}}}\\
&=(-1)^{r-1}\Big(\frac pq\Big)^r(1-w)
\sum_{h\ge1}\frac{wq^{h}w(1-q^{h-1})}{(1-wq^{h-1})^2(1-wq^h)}
\sum_{h<i_1<\dots <i_{r-1}}\frac{wq^{i_{1}}}{1-wq^{i_1}}\dots
\frac{wq^{i_{r-1}}}{1-wq^{i_{r-1}}}.
\end{align*}

Therefore
\begin{align*}
E^{(r)}_n&=(-1)^n[w^n](-1)^{r-1}\Big(\frac pq\Big)^r(1-w)^n
\sum_{h\ge1}\frac{wq^{h}w(1-q^{h-1})}{(1-wq^{h-1})^2(1-wq^h)}
\times\\&
\qquad\qquad\qquad\qquad\qquad\qquad
\times
\sum_{h<i_1<\dots <i_{r-1}}\frac{wq^{i_{1}}}{1-wq^{i_1}}\dots
\frac{wq^{i_{r-1}}}{1-wq^{i_{r-1}}}\\
&=(-1)^{r-1}\Big(\frac pq\Big)^r\sum_{k=r+1}^n\binom nk(-1)^k
\times\\&
\qquad\qquad\qquad
\times
[w^k]\sum_{1\le h<i_1<\dots <i_{r-1}}
\frac{wq^{h}w(1-q^{h-1})}{(1-wq^{h-1})^2(1-wq^h)}
\frac{wq^{i_{1}}}{1-wq^{i_1}}\dots
\frac{wq^{i_{r-1}}}{1-wq^{i_{r-1}}}.
\end{align*}

For the computation of the inner sum we need
\begin{multline*}
[w^k]
\sum_{h\ge1}\frac{(wq^{h})^{M}}{(1-wq^{h-1})^2(1-wq^h)}\\
=-\frac{q^{M+1}}{(1-q)^2}
\frac{1}{1-q^k}+\frac{q^M}{1-q}
\frac{1}{1-q^k}(k+1-M)+\frac{q^2}{(1-q)^2}
\frac1{Q^k-1}
\end{multline*}

and also (as a consequence)

\begin{align*}
[w^k]
\sum_{h\ge1}&\frac{(wq^{h})^{M}w(1-q^{h-1})}{(1-wq^{h-1})^2(1-wq^h)}\\
&=-\frac{q^{m+2}}{(1-q)^2}
\frac{1}{1-q^{k-1}}+\frac{q^{m+1}}{1-q}
\frac{1}{1-q^{k-1}}(k-m-1)+\frac{q^2}{(1-q)^2}
\frac1{Q^{k-1}-1}\\
&+\frac{q^{m+2}}{(1-q)^2}
\frac{1}{1-q^k}-\frac{q^{m+1}}{1-q}
\frac{1}{1-q^k}(k-m-1)-\frac{q^2}{(1-q)^2}
\frac1{Q^k-1}.
\end{align*}

The righthand side will be denoted by $\tau(k,m)$
and can alternatively be written as

\begin{align*}
\tau(k,m)&=-\frac{q^{m+2}}{p^2}
\frac{1}{Q^{k-1}-1}+\frac{q^{m+1}}{p}
\frac{1}{Q^{k-1}-1}(k-m-1)+\frac{q^2}{p^2}
\frac1{Q^{k-1}-1}\\*
&+\frac{q^{m+2}}{p^2}
\frac{1}{Q^k-1}-\frac{q^{m+1}}{p}
\frac{1}{Q^k-1}(k-m-1)-\frac{q}{p^2}
\frac1{Q^k-1}.\\
\end{align*}

Consequently

\begin{multline*}
[w^k]
\sum_{h\ge1}\frac{wq^{h}w(1-q^{h-1})}{(1-wq^{h-1})^2(1-wq^h)}
\sum_{h<i_1<\dots <i_{r-1}}\frac{wq^{i_{1}}}{1-wq^{i_1}}\dots
\frac{wq^{i_{r-1}}}{1-wq^{i_{r-1}}}\\
=\sum_{0<l_1<\dots<l_{r-1}< k-1}
\frac{1}{Q^{l_{1}}-1}\dots\frac{1}{Q^{l_{r-1}}-1}
\tau(k,l_{r-1}).
\end{multline*}

It is instructive to start with
the case $r=1$:

\begin{align*}
E^{(1)}_n&=\frac pq\sum_{k=2}^n\binom nk(-1)^k[w^k]
\sum_{h\ge1}\frac{wq^{h}w(1-q^{h-1})}{(1-wq^{h-1})^2(1-wq^h)}\\
&=\frac pq\sum_{k=2}^n\binom nk(-1)^k\Big[
\frac qp\frac{k-1}{Q^{k-1}-1}-
\frac qp\frac{k}{Q^{k}-1}\Big]\\
&=\sum_{k=2}^n\binom nk(-1)^k
\Big[\frac{k-1}{Q^{k-1}-1}-
\frac{k}{Q^{k}-1}\Big].\\
\end{align*}

When applying Rice's formula, it turns out that the
dominant singularity is at $z=1$.
This matches the intuition, since the position should be roughly
at  $\approx \mathit{const}\cdot n$.

Indeed, apart from fluctuations (coming from complex poles
at $z=1+2k\pi i/\log Q$, with $k\in \mathbb{Z}\setminus\{0\}$) we get
$E^{(1)}_n\sim n(\frac 1{L}-\frac qp)$.

Formally, for $q\to1$, we get  
$\frac n2$, as it should. 

Now, for the general case, define the function
\begin{equation*}
\psi_r(k)=\sum_{0<l_1<\dots<l_{r-1}< k-1}
\frac{1}{Q^{l_{1}}-1}\dots\frac{1}{Q^{l_{r-1}}-1}
\tau(k,l_{r-1}).
\end{equation*}

We need the value  $\psi_r(1)$, because of the result
(apart from fluctuations)
\begin{equation}\label{erw}   
E^{(r)}_n\sim
(-1)^{r-1}\Big(\frac pq\Big)^r \psi_r(1)\, n.
\end{equation}

We want to study the constants in more detail. 
Assume $r=2$. Since $\tau(z,z-1)=0$, 
the adaequate extension is

\begin{equation*}
\psi_2(z)=
\sum_{l\ge1}\frac1{Q^l-1}\tau(z,l)-
\sum_{l\ge1}\frac1{Q^{l+z-1}-1}\tau(z,l+z-1).
\end{equation*}

An elementary computation yields
\begin{align*}
\tau(z,l)&\sim
\frac{-(l-1)q^{l+2}+lq^{l+1}-q^2}{Lp^2}\frac{1}{z-1}\\*&
\qquad+
\frac{q^{l+1}}{pL}+
\frac{1}{2p^3}\left[-3q^2+q^3+lq^{l+1}+q^{l+2}-(l-1)q^{l+3}\right]+
\dots,\\
\tau(z,z+l-1)&\sim
\frac{-(l-1)q^{l+2}+lq^{l+1}-q^2}{Lp^2}\frac{1}{z-1}\\*&
\qquad+
\frac{1}{2p^3}\left[-3q^2+q^3+3lq^{l+1}-(4l-3)q^{l+2}+lq^{l+3}
\right],\\
\frac1{Q^{l+z-1}-1}&\sim\frac1{Q^{l}-1}+\frac{LQ^l}{(Q^l-1)^2}(z-1)+
\dots\; .
\end{align*}

Hence
\begin{align*}
\psi_2(1)&=\frac{1}{pL}\sum_{l\ge1}\frac{q^{l+1}}{Q^l-1}
+\frac{1}{p^3}\sum_{l\ge1}\frac1{Q^l-1}
\left[-
lq^{l+1}+(2l-1)q^{l+2}-(l-1)q^{l+3}
\right]\\
&+\frac{1}{p^2}\sum_{l\ge1}\frac{Q^l}{(Q^l-1)^2}
\left[(l-1)q^{l+2}-lq^{l+1}+q^2\right]\\
&=\frac{q}{pL}\sum_{l\ge1}\frac{1}{Q^l-1}-\frac{q^2}{p^2L}
+\frac{q^2(q+1)}{p^3}-\frac{q}{p}\sum_{l\ge1}\frac{lQ^l}{(Q^l-1)^2}
\end{align*}
and thus

\begin{align*}
-\frac{p^2}{q^2}\psi_2(1)
&=-\frac{p}{qL}\sum_{l\ge1}\frac{1}{Q^l-1}
+\frac pq\sum_{l\ge1}\frac{lQ^l}
{(Q^l-1)^2}+\frac1L-\frac 1p-\frac qp.\\
\end{align*}

It is our aim to show now that for $q\to1$  this 
quantity approaches $\frac{1}{4}$. 

For this, we write
 $q=e^{-t}$ and study the behaviour at $t\to0$. 
The {\sl Mellin transform\/} allows as to write

\begin{align*}
\sum_{l\ge1}\frac{1}{Q^l-1}
&=\frac{1}{2\pi i}\int_{3-i\infty}^{3+i\infty}\zeta^2(s)
\Gamma(s)t^{-s}ds,\\
\sum_{l\ge1}\frac{lQ^l}
{(Q^l-1)^2}&=\frac{1}{2\pi i}\int_{3-i\infty}^{3+i\infty}\zeta^2(s-1)
\Gamma(s)t^{-s}ds.\\
\end{align*}

Shifting the line of integration and collecting residues
gives us the following asymptotic expansions
(see \cite{FlGoDu95} for a friendly explanation of the
subject):
\begin{align*}
\sum_{l\ge1}\frac{1}{Q^l-1}
&=-\frac{\log t}{t}+\frac{\gamma}{t}+\frac{1}{4}+
\O(t),\\*
\sum_{l\ge1}\frac{lQ^l}
{(Q^l-1)^2}&=-\frac{\log t}{t^2}+\frac{1+\gamma}{t^2}
+\O(1).\\
\end{align*}

Using that, we find
\begin{align*}
-\frac{p^2}{q^2}\psi_2(1)
&=\frac{1}{4}-\frac{t}{9}+\O(t^2)=
\frac{1}{4}+\frac{1}{9}\log q+\O(\log^2  q).
\end{align*}

\bigskip

In the general case, the suitable extension is
\begin{equation*}   
\psi_r(z)=\sum_{l\ge1}g_{r-2}(l)\frac{\tau(z,l)}{Q^{l}-1}-
\sum_{l\ge1}g_{r-2}(l+z-1)\frac{\tau(z,l+z-1)}{Q^{l+z-1}-1},
\end{equation*}
with
\begin{equation*}   
g_s(l)=\sum_{1\le i_1<\dots<i_s<l}\frac{1}{Q^{i_1}-1}\dots
\frac{1}{Q^{i_s}-1},
\end{equation*}
and the suitable extension of it, which is, at least
theoretically, no problem, since we can, as already
explained in \cite{KnPr00}, use powersum symmetric
functions (see \cite{Stanley99}) or proceed as in the last section.

For instance, 
\begin{equation*}   
g_1(z)=\sum_{i\ge1}\frac{1}{Q^{i}-1}-
\sum_{i\ge1}\frac{1}{Q^{i+z-1}-1}
\end{equation*}
and
\begin{equation*}   
g_1(l+z-1)=g_1(l)+L
\sum_{i\ge l}\frac{Q^i}{(Q^{i}-1)^2}(z-1)+\dots,
\end{equation*}
 as $z\to1$. 

Using this, one could write an explicit formula for $\psi_3(1)$
which is already extremely messy.

Since already the quantity $\psi_2(1)$ cannot be
simplified, there is no hope for closed forms here. 

\bigskip

It is possible to make the translation to permutations
(the limit $q\to1$) at an earlier stage.

Note that
\begin{equation*}   
\lim_{q\to1}(1-q)\tau(k,m)=\frac{1}{2}-
\frac{(m+1)m}{2k}+\frac{m(m-1)}{2(k-1)}.
\end{equation*} 

Hence the constant in (\ref{erw}) translates into
\begin{equation*}   
\la_r:=(-1)^{r-1}\sum_{0<\l_1<\dots<l_{r-1}\le k-1}\frac{1}{l_1\dots l_{r-1}}
\bigg[\frac{1}{2}-
\frac{(l_{r-1}+1)l_{r-1}}{2k}+\frac{l_{r-1}(l_{r-1}-1)}{2(k-1)}\bigg]
\bigg|_{k=1}.
\end{equation*}

Observe that
\begin{equation*}   
\sum_{0<\l_1<\dots<l_{r-2}<l}\frac{1}{l_1\dots l_{r-2}}=
[x^{r-2}]\prod_{i=1}^{l-1}\left(1+\frac{x}{i}\right)=
[x^{r-2}]\frac{\Gamma(l+x)}{\Gamma(1+x)\Gamma(l)}.
\end{equation*}

Hence

\begin{align*}
\la_r&=-\frac{1}{2}[x^{r-2}]\sum_{1\le l \le k-1}
\frac{\Gamma(l-x)}{\Gamma(1-x)\Gamma(l)}
\left[\frac{1}{l}-\frac{l+1}{k}+\frac{l-1}{k-1}\right]\bigg|_{k=1}\\
&=[x^{r-2}]\frac{1}{2x}\left[
\frac{2\Gamma(k+1-x)}{\Gamma(k+1)\Gamma(3-x)}-1\right]\bigg|_{k=1}\\
&=[x^{r-2}]\frac{1}{4}\frac{1}{1-\frac x2}=2^{-r},
\end{align*}
as it should be.

In this permutation case we can even do better and compute
the next term (of order 1), because of
\begin{equation*}   
E^{(r)}_n\sim
(-1)^{r-1}\Big(\frac pq\Big)^r \Big(\psi_r(1)\, n-\psi_r(0)\Big).
\end{equation*}

All we have to do is to compute the following quantity:
\begin{align*}
\mu_r&=(-1)^{r-1}\sum_{0<\l_1<\dots<l_{r-1}\le k-1}\frac{1}{l_1\dots l_{r-1}}
\bigg[\frac{1}{2}-
\frac{(l_{r-1}+1)l_{r-1}}{2k}+\frac{l_{r-1}(l_{r-1}-1)}{2(k-1)}\bigg]
\bigg|_{k=0}\\
&=[x^{r-2}]\left(\frac{1}{1-x}-\frac{1}{4}\frac{1}
{1-\frac{x}{2}}\right)=1-2^{-r}.
\end{align*}

Since there are no more poles to the left and a technical
growth condition (see \cite{FlSe95}) is satisfied, there are
 no more terms, and we have the exact formula for the average
position of the $r$th left--to--right maximum, counted from the
right (remember that we counted one less than the position!):
\begin{equation*}   
E^{(r)}_n+1=\frac{n+1}{2^r}.
\end{equation*}

\begin{theorem}
The average position of  the      $ r$--th left--to--right
maximum from the right is given by
\begin{equation*}   
E^{(r)}_n=(-1)^r\Big(\frac pq\Big)^{r-1}\psi_r(1)\,n+
\varpi_r(\log_Qn)+o(n),
\end{equation*}
where the constants are defined in the text, and 
$\varpi_r(x)$ is a small periodic function, originating from the complex
poles $1+2k\pi i/\log Q$.

In the limit $q\to1$ (permutations), this simplifies to
 \begin{equation*}   
E^{(r)}_n=\frac n{2^r}+o(n).
\end{equation*}

\end{theorem}

 \bibliographystyle{plain}


\begin{thebibliography}{1}

\bibitem{Andrews76}
G.~Andrews.
\newblock {\em The Theory of Partitions}, volume~2 of {\em Encyclopedia of
  Mathematics and its Applications}.
\newblock Addison--Wesley, 1976.

\bibitem{FlGoDu95}
P.~Flajolet, X.~Gourdon, and P.~Dumas.
\newblock Mellin transforms and asymptotics: {H}armonic sums.
\newblock {\em Theoretical Computer Science}, 144:3--58, 1995.

\bibitem{FlSe95}
P.~Flajolet and R.~Sedgewick.
\newblock Mellin transforms and asymptotics: Finite differences and {R}ice's
  integrals.
\newblock {\em Theoretical Computer Science}, 144:101--124, 1995.

\bibitem{FlSe05}
P.~Flajolet and R.~Sedgewick.
\newblock {\em Analytic Combinatorics}.
\newblock Cambridge University Press, Cambridge, 2005.

\bibitem{KnPr00}
A.~Knopfmacher and H.~Prodinger.
\newblock Combinatorics of geometrically distributed random variables: Value
  and position of the $r$th left-to-right maximum.
\newblock {\em Discrete Mathematics}, 2000.

\bibitem{Prodinger96c}
H.~Prodinger.
\newblock Combinatorics of geometrically distributed random variables:
  Left-to-right maxima.
\newblock {\em Discrete Mathematics}, 153:253--270, 1996.

\bibitem{Stanley99}
R.~Stanley.
\newblock {\em Enumerative combinatorics. {V}ol. 2}.
\newblock Cambridge University Press, Cambridge, 1999.

\bibitem{SzRe90}
W.~Szpankowski and V.~Rego.
\newblock Yet another application of a binomial recurrence. {O}rder statistics.
\newblock {\em Computing}, 43(4):401--410, 1990.

\end{thebibliography}
\end{document}